\documentclass[12pt,leqno,a4paper]{article}
\usepackage{amsfonts}
\usepackage{amsmath, amssymb, amsthm, amsfonts}

\numberwithin{equation}{section}

\begin{document}
\title{On the  mean square of the error term for the two-dimensional divisor problem (I)}
\author{ Wenguang  Zhai$^1$,  Xiaodong Cao$^2$\\ \\
1   Shandong Normal University, Jinan, P. R. China\\
2 Beijing Institute of Petro-Chemical Technology, Beijing, P. R.
China}
\date{}

 \footnotetext[0]{This
work is supported by National Natural Science Foundation of
China(Grant No. 10771127) and National Natural Science Foundation of
Shandong Province(Grant No. Y2006A31).} \maketitle

{\bf Abstract.} Suppose $a$ and $b$ are two fixed positive integers
such that $(a,b)=1.$ In this paper we shall establish an asymptotic
formula for  the mean square of the error term $\Delta_{a,b}(x)$ of
the general two-dimensional divisor problem.

\section{\bf Introduction }

Suppose $1\leq a\leq b$ are two fixed integers. Without loss of
generality, we suppose  $(a,b)=1.$ Define
$d_{a,b}(n):=\sum_{n=h^ar^b}1.$
 The general two-dimensional divisor problem is to study the
error term
\begin{eqnarray*}
\Delta_{a,b}(x):&&= \sum_{n\leq x}
d_{a,b}(n)-\zeta(b/a)x^{1/a}-\zeta(a/b)x^{1/b},
\end{eqnarray*}
if $a\not= b.$ If $a=b,$ then the appropriate limit is to be taken
in the above sum. This problem attracts the interests of many
authors.

 When $a=b=1,$ $\Delta_{1,1}(x)$ is the error term of the well-known
Dirichlet divisor problem. Dirichlet first proved that
 $\Delta_{1,1}(x)=O(x^{1/2}).$
 The exponent $1/2$ was improved by many authors.
The latest result reads(see Huxley\cite{Hu})
\begin{equation}
\Delta_{1,1}(x)\ll  x^{131/416}(\log x)^{26947/8320 }.
\end{equation}
For the lower bounds, the best results read
\begin{equation}
\Delta_{1,1}(x)=\Omega_+(x^{1/4}(\log x)^{1/4}(\log\log x)^{(3+\log
4)/4} exp(-c\sqrt {\log\log\log x})) (c>0)
\end{equation}
 and
\begin{equation}
\Delta_{1,1}(x)=\Omega_{-}(x^{1/4}exp(c^{'}(\log\log
x)^{1/4}(\log\log\log x)^{-3/4})) \hspace{3mm}(c^{'}>0),
\end{equation}
which are due to Hafner\cite{Ha1} and  Corr\'adi and
K\'atai\cite{CK}, respectively.

When $a\not= b,$ Richert\cite{Ri} proved that
\begin{eqnarray}
\Delta_{a,b}(x)\ll\left\{\begin{array}{ll}
x^{\frac{2}{3(a+b)}},&\mbox{ $b\leq 2a,$}\\
x^{\frac{2}{5a+2b}},&\mbox{ $b\geq 2a.$}
\end{array}\right.
\end{eqnarray}
Better upper estimates can be found in \cite{Kr1,Kr2, Ra,S}.
Hafner\cite{Ha2} proved that
\begin{equation}
\Delta_{a,b}(x)=\Omega_{+}(x^{1/2(a+b)}(\log x)^{b/2(a+b)}\log\log
x)
\end{equation}
and
\begin{equation}
\Delta_{a,b}(x)=\Omega_{-}(x^{1/2(a+b)} e^{U(x)} ),
\end{equation}
where
$$U(x)=B(\log\log x)^{b/2(a+b)}(\log\log\log x)^{b/2(a+b)-1}$$
for some positive constant $B>0.$

It is conjectured that the estimate
\begin{equation}
\Delta_{a,b}(x)=O(x^{1/2(a+b)+\varepsilon})
\end{equation}
holds for any $1\leq a\leq b, (a,b)=1. $ When $a=b=1, $ the
conjecture (1.7) is supported by the power moment results of
$\Delta_{1,1}(x).$  For the mean square of $\Delta_{1,1}(x),$
Cram\'er\cite{Cr} first proved the classical  result
\begin{equation}
\int_1^T\Delta_{1,1}^2(x)dx=\frac{(\zeta(3/2))^4}{6\pi^2
\zeta(3)}T^{3/2} +O(T^{5/4+\varepsilon} ).
\end{equation}
  The estimate $O(T^{5/4+\varepsilon} )$ in (1.8) was
improved to $O(T\log^5 T)$ in \cite{To} and $O(T\log^4 T)$ in
\cite{P}.   The higher-power moments of $\Delta_{1,1}(x)$ were
studied in\cite{H1,IS, Ts, Zw1, Zw2, Zw3}. When $a\not= b,$ by using
the theory of the Riemann zeta-function, Ivi\'c\cite{I4} proved that
\begin{equation}
\int_1^T\Delta_{a,b}^2(x)dx\left\{\begin{array}{ll}
\ll T^{1+1/(a+b)}\log^2 T, \\
=\Omega(T^{1+1/(a+b)}).
\end{array}\right.
\end{equation}
He also conjectured that the asymptotic formula
\begin{equation}
\int_1^T\Delta_{a,b}^2(x)dx=c_{a,b}T^{1+1/(a+b)}(1+o(1))
\end{equation}
 holds for some positive constant $c_{a,b}.$ We note that the
 $\Omega$-result in (1.9) is also contained in a very general result
 of the second-named author(see Theorem 5 of \cite{Cx}).

 When $a=b=1,$ the proofs of most
power moment results of $\Delta_{1,1}(x)$ mentioned above were
started from the well-known truncated Voronoi's formula(see, for
example, \cite{I1})
\begin{equation}
\Delta_{1,1}(x)= (\sqrt 2\pi)^{-1}x^{1/4}\sum_{n\leq N}\frac{d(n)}
{n^{3/4}}\cos(4\pi\sqrt{xn}-\frac{\pi}{4})+O(x^\varepsilon+x^{1/2+\varepsilon}N^{-1/2})
\end{equation}
for $1\ll N\ll x^A,$ where $A>0$ is any fixed constant. Note that
the infinite series $\sum_{n=1}^\infty d(n)
n^{-3/4}\cos(4\pi\sqrt{xn}-\pi/4)$ is conditionally convergent. When
$a\not= b,$ Kr\"{a}tzel\cite{Kr1} provided a series representation
of $\Delta_{a,b}(x).$ However, Kr\"{a}tzel's series converges only
when $b<3a/2,$ which is a significant restriction for many
applications.

The aim of this paper is to prove Ivi\'c's conjecture (1.10). More
precisely, we shall prove the  following Theorem.

{\bf Theorem .} Suppose $1\leq a< b$ are fixed integers for which
$(a,b)=1.$  Then we have
\begin{eqnarray}
\int_1^T\Delta_{a,b}^2(x)dx=c_{a,b}T^{\frac{1+a+b}{a+b}}+O(T^{\frac{1+a+b}{a+b}-\frac{a}{2b(a+b)(a+b-1)}
}\log^{7/2} T),
\end{eqnarray}
where
\begin{eqnarray*}
&&c_{a,b}:=\frac{a^{b/(a+b)}b^{a/(a+b)}}{2(a+b+1)\pi^2}\sum_{n=1}^\infty
g_{a,b}^2(n),\\
&&g_{a,b}(n):=\sum_{n=h^ar^b}h^{-\frac{a+2b}{2a+2b}}r^{-\frac{b+2a}{2a+2b}}.
\end{eqnarray*}

{\bf Remark 1.} It is easy to see that the function $g_{a,b}(n)$ is
symmetric for $a$ and $b$,
 namely   $g_{a,b}(n)=g_{b,a}(n).$ The convergence of the infinite series
$\sum_{n=1}^\infty g^2_{a,b}(n)$ will be proved in Section 4.

 {\bf Remark 2.} Our theorem  also holds for
$a=b=1.$ In this case we have $$g_{1,1}(n)=d(n)n^{-3/4},\
c_{1,1}=\frac{1}{6\pi^2}\sum_{n=1}^\infty
d^2(n)n^{-3/2}=\frac{(\zeta(3/2))^4}{6\pi^2 \zeta(3)},$$ where
$d(n)=d_{1,1}(n)$ is the Dirichlet divisor function. Hence our
Theorem provides a new proof of Cramer's  classical result  (1.8).

{\bf Notations.} ${\Bbb Z}$ denotes the set of all integers. For a
real number $u,$ $[u]$ denotes the integer part of $u,$ $\{u\}$
denotes the fractional part of $u,$ $\psi(u)=\{u\}-1/2,$ $\Vert
u\Vert$ denotes the distance from $u$ to the integer nearest to $u.$
$\mu(n)$ is the M\"obius function, $(m,n)$ denotes the greatest
common divisor of natural numbers $m$ and $n.$ $n\sim N$ means
$N<n\leq 2N.$ $\varepsilon$ always denotes a sufficiently small
positive constant. $SC(\Sigma)$ denotes the summation condition of
the sum $\Sigma$ when it is complicated. Finally, define
\begin{eqnarray*}
\sideset{}{^{\prime}}\sum_{\alpha\leq n\leq
\beta}f(n)=\left\{\begin{array}{ll}
\sum_{\alpha< n< \beta}f(n),&\mbox{$\alpha\notin {\Bbb Z}, \beta\notin {\Bbb Z},$}\\
f(\alpha)/2+\sum_{\alpha< n< \beta}f(n),&\mbox{$\alpha\in {\Bbb Z}, \beta\notin {\Bbb Z},$}\\
f(\alpha)/2+\sum_{\alpha< n< \beta}f(n)+f(\beta)/2,&\mbox{$\alpha\in
{\Bbb Z}, \beta\in {\Bbb Z}.$}
\end{array}\right.
\end{eqnarray*}

\section{\bf Two preliminary Lemmas}

In order to prove our theorem, we need the following two Lemmas.
Lemma 2.1 is well-known; see for example, Heath-Brown\cite{H2}.
Lemma 2.2 is Theorem 2.2 of Min\cite{Mi},  see also Lemma 6 of
Chapter 1 in \cite{V}. A weaker version of Lemma 2.2 can be found in
\cite{KN}, which also suffices for our proof.

{\bf Lemma 2.1.}   Let $H\geq 2$ be any real number. Then
$$\psi(u)=-\sum_{1\leq |h|\leq H}\frac{e(hu)}{2\pi ih}+O\left(\min(1,\frac{1}{H\Vert u\Vert}) \right).$$

{\bf Lemma 2.2.} Suppose  $A_1,\cdots, A_5$ are absolute positive
constants, $f(x)$ and $ g(x)$ are algebraic functions in $[a,b]$ and
\begin{eqnarray*}
&&\frac{A_1}{R}\leq |f^{''}(x)|\leq\frac{A_2}{R},\  \ \ |f^{'''}(x)|\leq\frac{A_3}{RU},\ \ U\geq 1, \\
&&|g(x)|\leq A_4G,\ \ \  |g^{'}(x)|\leq A_5GU_1^{-1},\ \ U_1\geq 1,
\end{eqnarray*}
$[\alpha,\beta]$ is the image of  $[a,b]$ under the mapping
$y=f^{'}(x)$, then
\begin{eqnarray*}
\sum_{a<n\leq b}g(n)e(f(n))&=&e^{\pi i/4}\sum_{\alpha\leq u\leq
\beta}b_u
\frac{g(n_u)}{\sqrt{f^{''}(n_u)}}e\left(f(n_u)-un_u \right)\\
&&+O\left(G\log(\beta-\alpha+2)+G(b-a+R)(U^{-1}+U_1^{-1})\right)\\
&&+ O\left(G\min\left[\sqrt R,
\max\left(\frac{1}{<\alpha>},\frac{1}{<\beta>}\right)\right]\right),
\end{eqnarray*}
where  $n_u$  is the solution of $f^{'}(n)=u$,
\begin{eqnarray*}
<t>=\left\{\begin{array}{ll}
\Vert t\Vert,&\mbox{if $t$ not an integer,}\\
\beta-\alpha,& \mbox{if $t$ an integer,}
\end{array}\right.
\end{eqnarray*}
\begin{eqnarray*}
b_u=\left\{\begin{array}{ll}
1 ,&\mbox{if $\alpha<u<\beta$, or $\alpha, \beta$ not integers ,}\\
1/2,& \mbox{if $\alpha$ or $\beta$ are integers,}\\
\end{array}\right.
\end{eqnarray*}
\begin{eqnarray*}
\sqrt{f^{\prime\prime}}=\left\{\begin{array}{ll}
 \sqrt{f^{\prime\prime}},&\mbox{if $ f^{\prime\prime}>0,$}\\
 i\sqrt{|f^{\prime\prime}|},& \mbox{if $ f^{\prime\prime}<0. $}
\end{array}\right.
\end{eqnarray*}

\section{\bf A Voronoi type formula of $\Delta(a,b; x)$}

It suffices for us to evaluate the integral
$\int_T^{2T}\Delta_{a,b}^2(x)dx ,$ where $T\geq 10$ is a large
parameter.

It is well-known that(see for example, Ivi\'c\cite{I3}, eq.(14.46))
\begin{eqnarray}
\Delta_{a,b}(x)=f(a,b; x)+f(b,a; x)+O(1),
\end{eqnarray}
 where $$f(a,b; x):=-\sum_{m\leq
x^{1/(a+b)}}\psi\left(\frac{x^{1/a}}{m^{b/a}}\right).$$

Suppose   $T\leq x\leq 2T$,  $H$ is a parameter such that
$T^\varepsilon\ll H\ll T^{100(a+b)}.$  By Lemma 2.1 we have
\begin{eqnarray}
f(a,b; x)&&=R_1(a,b; x)+R_2(a,b; x),\\
R_1(a,b; x):&&=\frac{1}{2\pi i}\sum_{1\leq |h|\leq
H}\frac{1}{h}\sum_{m\leq
x^{1/(a+b)}}e\left(\frac{hx^{1/a}}{m^{b/a}}\right),\nonumber\\
R_2(a,b; x):&&=O\left(\sum_{m\leq
x^{1/(a+b)}}\min\left(1,\frac{1}{H\Vert \frac{x^{1/a}}{m^{b/a}}
\Vert}\right) \right).\nonumber
\end{eqnarray}

Define
\begin{eqnarray*}
&&c:=(2ab)^{ab}, J:=[ ({\cal L}/(a+b)-\log {\cal L})\log^{-1} c] ,
{\cal L}:=\log T,\\&& m_j:=x^{1/(a+b)}c^{-j}\ (j\geq
0).\end{eqnarray*}

It is easy to see that
$$c^J\asymp T^{1/(a+b)}{\cal L}^{-1}.$$

We have
\begin{eqnarray}
R_1(a,b; x)&&=\frac{1}{2\pi i}\sum_{1\leq |h|\leq
H}\frac{1}{h}\sum_{j=0}^J\sum_{m_{j+1}<m\leq
m_j}e\left(\frac{hx^{1/a}}{m^{b/a}}\right)+O({\cal L}^2)\\
&&=\frac{1}{2\pi i}\sum_{-H\leq h\leq
-1}\frac{1}{h}\sum_{j=0}^J\sum_{m_{j+1}<m\leq
m_j}e\left(\frac{hx^{1/a}}{m^{b/a}}\right)\nonumber\\
&&\ \ \ +\frac{1}{2\pi i}\sum_{1\leq h\leq
H}\frac{1}{h}\sum_{j=0}^J\sum_{m_{j+1}<m\leq
m_j}e\left(\frac{hx^{1/a}}{m^{b/a}}\right)+O({\cal L}^2)\nonumber\\
&&=-\frac{1}{2\pi i}\sum_{1\leq h\leq
H}\frac{1}{h}\sum_{j=0}^J\sum_{m_{j+1}<m\leq
m_j}e\left(-\frac{hx^{1/a}}{m^{b/a}}\right)\nonumber\\
&&\ \ \ +\frac{1}{2\pi i}\sum_{1\leq h\leq
H}\frac{1}{h}\sum_{j=0}^J\sum_{m_{j+1}<m\leq
m_j}e\left(\frac{hx^{1/a}}{m^{b/a}}\right)+O({\cal L}^2)\nonumber\\
&&=-\frac{\Sigma_1}{2\pi i}+\frac{\overline{\Sigma_1}}{2\pi
i}+O({\cal L}^2),\nonumber
\end{eqnarray}
say, where
$$\Sigma_1= \sum_{1\leq h\leq
H}\frac{1}{h}\sum_{j=0}^J\ \ \sum_{m_{j+1}<m\leq
m_j}e\left(-\frac{hx^{1/a}}{m^{b/a}}\right).$$

Let
$$S_{h,j}(x):=\sum_{m_{j+1}<m\leq
m_j}e\left(-\frac{hx^{1/a}}{m^{b/a}}\right).$$

Define
\begin{eqnarray*}
&&c_1(a,b):=a^{b/2(a+b)}b^{a/2(a+b)}(a+b)^{-1/2},\\
&&c_2(a,b):=a^{b/(a+b)}b^{-b/(a+b)}+b^{a/(a+b)}a^{-a/(a+b)}.
\end{eqnarray*}
It is easy to check that
\begin{equation}
c_1(a,b)=c_1(b,a),\ \  c_2(a,b)=c_2(b,a).
\end{equation}

By Lemma 2.2 we get
\begin{eqnarray}
S_{h,j}(x)&&=c_1(a,b)x^{\frac{1}{2(a+b)}}\sideset{}{^{\prime}}\sum_{n_{j,h}(a,b)\leq
r\leq
n_{j+1,h}(a,b)}h^{\frac{a}{2(a+b)}}r^{-\frac{2a+b}{2(a+b)}}\\
&&\ \ \ \ \ \
 \times
 e\left(-c_2(a,b)x^{\frac{1}{a+b}}(h^ar^b)^{\frac{1}{a+b}}-\frac{1}{8}\right)\nonumber
 +O({\cal L}),
\end{eqnarray}
where
$$n_{j,h}(a,b):=\frac{b}{a}h(2ab)^{(a+b)bj}.$$

Inserting (3.5) into $\Sigma_1$ we get
\begin{eqnarray}
\Sigma_1&&=c_1(a,b)x^{\frac{1}{2(a+b)}}\sum_{1\leq h\leq H}\ \ \
\sideset{}{^{\prime}}\sum_{\frac{bh}{a}\leq r\leq
n_{J+1,h}(a,b)}h^{-\frac{a+2b}{2(a+b)}}r^{-\frac{2a+b}{2(a+b)}}\\
&&\ \ \ \ \ \ \ \ \ \
 \times
 e\left(-c_2(a,b)x^{\frac{1}{a+b}}(h^ar^b)^{\frac{1}{a+b}}-\frac{1}{8}\right)\nonumber
 +O({\cal L}^2).
\end{eqnarray}

From (3.3) and (3.6) we get
\begin{eqnarray}
R_1(a,b; x)&&=R_1^{*}(a,b; x)+O({\cal L}^2),
\end{eqnarray}
where
\begin{eqnarray*}
\ \ \ R_1^{*}(a,b; x):&&=
\frac{c_1(a,b)}{\pi}x^{\frac{1}{2(a+b)}}\sum_{1\leq h\leq H}\ \ \
\sideset{}{^{\prime}}\sum_{\frac{bh}{a}\leq r\leq
n_{J+1,h}(a,b)}h^{-\frac{a+2b}{2(a+b)}}r^{-\frac{2a+b}{2(a+b)}}\\
&&\ \ \ \ \ \ \ \ \ \
 \times
 \cos\left(2\pi
 c_2(a,b)x^{\frac{1}{a+b}}(h^ar^b)^{\frac{1}{a+b}}-\frac{\pi}{4}\right)   .
\end{eqnarray*}

Define
\begin{eqnarray*}
g(a,b;n, H, J): &&=\sideset{}{^{\prime}}\sum_{\stackrel{n=h^ar^b,
1\leq h\leq H}{bh/a\leq r\leq
n_{J+1,h}(a,b)}}h^{-\frac{a+2b}{2(a+b)}}r^{-\frac{2a+b}{2(a+b)}},\\
g(a,b;n): &&=\sideset{}{^{\prime}}\sum_{\stackrel{n=h^ar^b
}{bh/a\leq r }}h^{-\frac{a+2b}{2(a+b)}}r^{-\frac{2a+b}{2(a+b)}}.
\end{eqnarray*}
It is easy to check that if    $h^ar^b\leq \min(H^{a+b},
T^{b/a}){\cal L}^{-b-b^2/a-1}, bh/a\leq r, $ then  it follows that
$h\leq H, r\leq n_{J+1,h}(a,b).$ Thus
\begin{equation}
g(a,b;n, H, J)=g(a,b;n),\ \ n\leq \min(H^{a+b}, T^{b/a}){\cal
L}^{-b-b^2/a-1}.
\end{equation}

So we have that
\begin{eqnarray*}
R_1^{*}(a,b; x)&&=
\frac{c_1(a,b)}{\pi}x^{\frac{1}{2(a+b)}}\sum_{1\leq n\leq
H^an_{J+1,H}(a,b)}g(a,b;n,H,J)\\&&\ \ \ \ \ \ \ \times\cos\left(2\pi
c_2(a,b)x^{\frac{1}{a+b}}n^{\frac{1}{a+b}}-\frac{\pi}{4}\right)
\end{eqnarray*}

Similarly we have
\begin{equation}
R_1(b,a;x)=R_1^{*}(b,a;x)+O({\cal L}^2),
\end{equation}
where
\begin{eqnarray*}
R_1^{*}(b,a; x)&&=
\frac{c_1(b,a)}{\pi}x^{\frac{1}{2(a+b)}}\sum_{1\leq n\leq
H^bn_{J+1,H}(b,a)}g^{*}(b,a;n,H,J)\\&&\ \ \ \ \ \ \
\times\cos\left(2\pi
c_2(b,a)x^{\frac{1}{a+b}}n^{\frac{1}{a+b}}-\frac{\pi}{4}\right),\\
 g(b,a;n, H, J):
&&=\sideset{}{^{\prime}}\sum_{\stackrel{n=h^br^a, 1\leq h\leq
H}{ah/b\leq r\leq
n_{J+1,h}(b,a)}}h^{-\frac{b+2a}{2(a+b)}}r^{-\frac{2b+a}{2(a+b)}},\\
g(b,a;n): &&=\sideset{}{^{\prime}}\sum_{\stackrel{n=h^br^a
}{ah/b\leq r }}h^{-\frac{b+2a}{2(a+b)}}r^{-\frac{2b+a}{2(a+b)}}.
\end{eqnarray*}
It is easy to check that if  that if $h^br^a\leq \min(H^{a+b},
T^{a/b}){\cal L}^{-a-a^2/b-1}, ah/b\leq r, $ then  $h\leq H, r\leq
n_{J+1,h}(b,a).$ Thus
\begin{equation}
g^{*}(b,a;n, H, J)=g(b,a;n),\ \ n\leq \min(H^{a+b}, T^{a/b}){\cal
L}^{-a-a^2/b-1}.
\end{equation}

Suppose $z$ is a parameter such that $T^\varepsilon\ll z\leq
\min(H^{a+b}, T^{a/b}){\cal L}^{-b-b^2/a-1}$ and define
\begin{eqnarray*}
R_{11}^{*}(a,b; x):&&=
\frac{c_1(a,b)}{\pi}x^{\frac{1}{2(a+b)}}\sum_{1\leq n\leq z
}g(a,b;n)\\&&\ \ \ \ \ \ \ \times\cos\left(2\pi
c_2(a,b)x^{\frac{1}{a+b}}n^{\frac{1}{a+b}}-\frac{\pi}{4}\right)\\
R_{12}^{*}(a,b; x):&&=R_{1}^{*}(a,b; x)-R_{11}^{*}(a,b; x),\\
R_{11}^{*}(b,a; x):&&=
\frac{c_1(b,a)}{\pi}x^{\frac{1}{2(a+b)}}\sum_{1\leq n\leq z
}g(b,a;n)\\&&\ \ \ \ \ \ \ \times\cos\left(2\pi
c_2(b,a)x^{\frac{1}{a+b}}n^{\frac{1}{a+b}}-\frac{\pi}{4}\right)\\
R_{12}^{*}(b,a; x):&&=R_{1}^{*}(b,a; x)-R_{11}^{*}(b,a; x).
\end{eqnarray*}

 Recalling (3.4) we have
\begin{eqnarray*}
R_{11}^{*}(a,b; x)+R_{11}^{*}(b,a; x)&&=
\frac{c_1(a,b)}{\pi}x^{\frac{1}{2(a+b)}}\sum_{1\leq n\leq z
}\left(g(a,b;n)+g(b,a;n)\right)\\&& \ \ \ \ \  \times\cos\left(2\pi
c_2(a,b)x^{\frac{1}{a+b}}n^{\frac{1}{a+b}}-\frac{\pi}{4}\right).
\end{eqnarray*}

From the definition of $g(a,b;n)$ and $g(b,a;n)$ we have
\begin{eqnarray*}
&&\ \ \ \ g(a,b;n)+g(b,a;n)\\
&&=\sideset{}{^{\prime}}\sum_{\stackrel{n=h^ar^b }{bh/a\leq r
}}h^{-\frac{a+2b}{2(a+b)}}r^{-\frac{2a+b}{2(a+b)}} +
\sideset{}{^{\prime}}\sum_{\stackrel{n=h^br^a }{ah/b\leq r
}}h^{-\frac{b+2a}{2(a+b)}}r^{-\frac{2b+a}{2(a+b)}}\\
&&=\sideset{}{^{\prime}}\sum_{\stackrel{n=h^ar^b }{bh/a\leq r
}}h^{-\frac{a+2b}{2(a+b)}}r^{-\frac{2a+b}{2(a+b)}} +
\sideset{}{^{\prime}}\sum_{\stackrel{n=r^bh^a }{ar/b\leq h
}}r^{-\frac{b+2a}{2(a+b)}}h^{-\frac{2b+a}{2(a+b)}}\\
&&=\sum_{ n=h^ar^b}h^{-\frac{a+2b}{2(a+b)}}r^{-\frac{2a+b}{2(a+b)}}=
g_{a,b}(n).
\end{eqnarray*}

Define
\begin{eqnarray*}
\Delta^{*}_{a,b}(x,z)&&:=
\frac{c_1(a,b)}{\pi}x^{\frac{1}{2(a+b)}}\sum_{1\leq n\leq z
}g_{a,b}(n)\\&&\ \ \ \ \ \ \ \times\cos\left(2\pi
c_2(a,b)x^{\frac{1}{a+b}}n^{\frac{1}{a+b}}-\frac{\pi}{4}\right).
\end{eqnarray*}

Combining the above estimates  we get
\begin{eqnarray}
\Delta_{a,b}(x)&&=\Delta^{*}_{a,b}(x,z)+E_{a,b}(x),\\
E_{a,b}(x):&&=R_{12}^{*}(a,b;x)+R_{12}^{*}(b,a;x)+R_2(a,b;x)\nonumber\\
&&\ \ \ \ \ +R_2(b,a;x)+O({\cal L}^2).\nonumber
\end{eqnarray}

The formula (3.11) can be viewed as a truncated Voronoi's formula.

\section{\bf On the series $\sum_{n=1}^\infty g_{a,b}^2(n)$}\

In this  section we shall prove that the infinite  series
$\sum_{n=1}^\infty g^2_{a,b}(n)$ is convergent. Without loss of
generality, we suppose $a<b.$ By the definition of $g_{a,b}(n)$ it
is easy to see that
\begin{eqnarray}
g_{a,b}(n)=n^{-\frac{2a+b}{(2a+2b)b}}\sum_{n=h^ar^b}h^{-\frac{b-a}{a}}=n^{-\frac{2a+b}{(2a+2b)b}}g_{a,b}^{*}(n),
\end{eqnarray}
say, where  $$g_{a,b}^{*}(n): =\sum_{n=h^ar^b}h^{-\frac{b-a}{a}}$$
is a multiplicative function. So we have
\begin{eqnarray}
\sum_{n=1}^\infty\frac{g_{a,b}^{*2}(n)}{n^s}=\prod_{p}\left(1+\sum_{\alpha=1}^\infty\frac{g_{a,b}^{*2}(p^\alpha)}{p^{\alpha
s}}\right)\ \ (\Re s>1).
\end{eqnarray}

We shall show that
\begin{eqnarray}
&&\ \ \ \ \sum_{n=1}^\infty\frac{g_{a,b}^{*2}(n)}{n^s}=
\zeta(bs)\zeta(as+\frac{2(b-a)}{a}) G_{a,b}(s),
\end{eqnarray}
where $G_{a,b}(s)$ is analytic for $\Re s>(3a-2b)/a^2.$

We consider two cases  .

{\bf Case 1.} $a=1$

If $1\leq \alpha\leq b-1,$ then $p^\alpha=h^ar^b$ implies
$h=p^\alpha,r=1.$ So
$$g_{a,b}^{*}(p^\alpha)=p^{-\alpha (b-1)}.$$

Suppose  $\alpha=mb$ with $m\geq 1.$ Then $p^\alpha=h^ar^b$ implies
  $h=p^{jb}, r=p^{(m-j)b}, j=0,1,2,\cdots, m$. So
$$g_{a,b}^{*}(p^\alpha) =\sum_{j=0}^mp^{-jb(b-1)}.$$

Suppose $\alpha=mb+u$ with $m\geq 1$ and $1\leq u\leq b-1.$ Then
$p^\alpha=h^ar^b$ implies  $  h=p^{jb+u}, r=p^{(m-j)b},
j=0,1,2,\cdots, m. $ So
$$g_{a,b}^{*}(p^\alpha) =\sum_{j=0}^m p^{-(jb+u)(b-1)}=p^{-u(b-1)}+p^{-u(b-1)}\sum_{j=1}^mp^{-jb(b-1)}.$$

From the above we get($p\geq 2,\Re s=\sigma>1$)
\begin{eqnarray*}
&& \ \ \ \ \ \ 1+\sum_{\alpha=1}^\infty\frac{
g_{a,b}^{*2}(p^\alpha)}{p^{\alpha
s}} \\
&&= 1+\sum_{\alpha=1}^{b-1}\frac{p^{-2\alpha(b-1)}}{p^{\alpha s}}
+\sum_{m=1}^\infty\frac{(\sum_{j=0}^mp^{-jb(b-1)})^2
}{p^{mbs}}+\sum_{m=1}^\infty\sum_{u=1}^{b-1} \frac{(\sum_{j=0}^m
p^{-(jb+u)(b-1)})^2}{p^{(mb+u)s}} \nonumber\\
 &&=1+p^{-bs}+O(p^{-\sigma-2(b-1)}).
\end{eqnarray*}
Inserting this formula into (4.2) we get that (4.3) is true for
$a=1. $

{\bf Case 2.} $a\geq 2$

Let $k=[b/a].$ Then $ (a,b)=1$ implies that $ka<b<ka+a.$

Consider the equation
\begin{equation}
ua+vb=\alpha ,\ \ u,v\in {\Bbb Z}.
\end{equation}

Let $v(\alpha;a,b)$ denote the number of non-negative solutions of
the equation (4.4). If $v(\alpha;a,b)=0,$ then $
g_{a,b}^{*}(p^\alpha)=0.$ Especially , it is easy to see that if
  $\alpha\leq b-1, a\not| \alpha,$ then $v(\alpha;a,b)=0.$ Hence
$$ g_{a,b}^{*}(p^\alpha)=0,\ \ (\alpha\leq b-1, a\not| \alpha).$$

If $\alpha=ja,(j=1,2,\cdots, k),$ then the equation (4.4) has only
one solution $(u,v)=(j,0).$ So $p^\alpha=h^ar^b$ implies that
$h=p^j,r=1,$ and hence
$$ g_{a,b}^{*}(p^\alpha)=p^{-\frac{j(b-a)}{a}}.$$

Later suppose $\alpha\geq b$ and $v(\alpha;a,b)>0.$ Let $(u_0,v_0)$
denote a special non-negative solution of (4.4) such that $0\leq
u_0<b.$ Then all non-negative solutions of (4.4) are
$(u,v)=(u_0+bt,v_0-at), 0\leq t\leq v(\alpha;a,b).$ Hence we have
$$g_{a,b}^{*}(p^\alpha) =\sum_{j=0}^{v(\alpha;a,b)-1}p^{-\frac{(u_0+bj)(b-a)}{a}} .$$
Especially when $b|\alpha,$ $(0,\alpha/b)$ is a non-negative
solution of (4.4). Thus
$$g_{a,b}^{*}(p^\alpha) =\sum_{j=0}^{v(\alpha;a,b)-1}p^{-\frac{bj(b-a)}{a}} ,\ \ \alpha\equiv 0(mod\ b).$$

From the above  we get
\begin{eqnarray*}
&& \ \ \ \ \ 1+\sum_{\alpha=1}^\infty\frac{
g_{a,b}^{*2}(p^\alpha)}{p^{\alpha
s}} \\
&&= 1+\sum_{j=1}^{k}\frac{p^{-2j(b-a)/a}}{p^{jas}} +\sum_{\alpha>b,
b|\alpha}\frac{ g_{a,b}^{*2}(p^\alpha)}{p^{\alpha s}}
+\sum_{\alpha>b, b \not|\alpha}\frac{
g_{a,b}^{*2}(p^\alpha)}{p^{\alpha s}}\\
&&=1+p^{-bs}+O(p^{-a\sigma-\frac{2(b-a)}{a}}).
\end{eqnarray*}

Inserting this formula into (4.2) we get that (4.3) is true in this
case.

\bigskip

It is easy to check that $\zeta(as+\frac{2b-2a}{a})$ has a simple
pole at $s=(3a-2b)/a^2. $ Since $b>a\geq 1,$ it is easy to see that
$(3a-2b)/a^2<1/b.$ From (4.3)  and Perron's formula we get that for
any $X>2, $
\begin{eqnarray*}
\sum_{X<n\leq 2X} g_{a,b}^{*2}(n)\ll X^{1/b}.
\end{eqnarray*}
Hence we get
\begin{eqnarray}
\sum_{X<n\leq 2X} g_{a,b}^2(n)\ll X^{\frac
1b-\frac{2a+b}{(a+b)b}}\ll X^{-\frac{a}{(a+b)b}},
\end{eqnarray}
which implies   the convergence of the infinite series
$\sum_{n=1}^\infty g_{a,b}^2(n)$  . From (4.5) we also get
\begin{eqnarray}
\sum_{n\leq X} g_{a,b}^2(n)= \sum_{n=1}^\infty
g^2_{a,b}(n)+O(X^{-\frac{a}{(a+b)b}})
\end{eqnarray}
and
\begin{eqnarray}
\sum_{n>X} g^2_{a,b}(n)= O(X^{-\frac{a}{(a+b)b}}).
\end{eqnarray}

\section{\bf On a special sum}

Suppose $a, b$ are fixed natural integers, $T$ is a large parameter
. Define
\begin{eqnarray*}
S_{a,b}(T):&&=\sum_{2}(h_1h_2)^{-\frac{2b+a}{2a+2b}}(r_1r_2)^{-\frac{2a+b}{2a+2b}}\\&&
 \ \ \ \ \times
\min\left(T^{\frac{1}{a+b}},\frac{1}{\left|h_1^{\frac{a}{a+b}}r_1^{\frac{b}{a+b}}-h_2^{\frac{a}{a+b}}r_2^{\frac{b}{a+b}}\right|}\right),
\end{eqnarray*}
where
\begin{eqnarray*}
SC(\Sigma_2): &&h_1^ar_1^b\leq T^{100(a+b)}, h_2^ar_2^b\leq T^{100(a+b)},\\
&& 0<|h_1^{\frac{a}{a+b}}r_1^{\frac{b}{a+b}}-
h_2^{\frac{a}{a+b}}r_2^{\frac{b}{a+b}}|<\frac{1}{10}
h_1^{\frac{a}{2a+2b}}r_1^{\frac{b}{2a+2b}}
 h_2^{\frac{a}{2a+2b}}r_2^{\frac{b}{2a+2b}}.
\end{eqnarray*}

In this section we shall estimate the sum $S_{a,b}(T)$, which is
very important in our proof.

\subsection{\bf On a Diophantine inequality}\

Suppose $\alpha$ and $\beta$ are fixed non-zero real numbers,
$H_1\geq 1, H_2\geq 1, R_1\geq 1, R_2\geq 1$ are large real numbers,
$  \delta>0.$ Let ${\cal A}(H_1, H_2, R_1, R_2;\delta)$ denote the
number of the solutions of the inequality
\begin{equation}
\left|  h_1^\alpha r_1^\beta- h_2^\alpha r_2^\beta\right|\leq
\delta,h_1\sim H_1, h_2\sim H_2, r_1\sim R_1, r_2\sim R_2.
\end{equation}

{\bf Lemma 5.1.} We have
\begin{align*}
{\cal A}(H_1, H_2, R_1, R_2;\delta)  &\ll\delta
(H_1H_2)^{1-\alpha/2}(R_1R_2)^{1-\beta/2}\\
&\ \ +(H_1H_2R_1R_2)^{1/2}(\log 2H_1H_2R_1R_2)^2,
\end{align*} where
the implied constant depends only on $\alpha, \beta.$

{\bf Remark.} When $H_1=H_2, R_1=R_2,$   Lemma 5.1 reduces to Lemma
1 of Fouvry and Iwaniec\cite{FI}. Here Lemma 5.1 is more general.

\begin{proof}
We follow the proof of Lemma 1 of Fouvry and Iwaniec \cite{FI}.
Suppose $u$ and $v$ are two positive integers and let ${\cal
A}_{u,v}(H_1,H_2,R_1, R_2;\delta)$ denote the number of  solutions
of the inequality (5.1) with $(r_1,r_2)=u,(h_1,h_2)=v.$ Set
$r_j=m_ju,h_j=l_jv (j=1,2),$ then $(m_1,m_2,l_1,l_2)$ satisfies
\begin{eqnarray}
&&\left|\frac{m_1^\beta}{m_2^\beta}-\frac{l_2^\alpha}{l_1^\alpha}\right|\leq
c(\alpha)c(\beta)\delta R_2^{-\beta}H_1^{-\alpha},\\
&&\left|\frac{m_2^\beta}{m_1^\beta}-\frac{l_1^\alpha}{l_2^\alpha}\right|\leq
c(\alpha)c(\beta)\delta R_1^{-\beta}H_2^{-\alpha},
\end{eqnarray}
where
\begin{eqnarray*}
c(\gamma)=\left\{\begin{array}{ll}
1,&\mbox{if  $\gamma>0,$}\\
2^{-\gamma},& \mbox{if $\gamma<0.$}
\end{array}\right.
\end{eqnarray*}

It is easy to see that $\frac{m_1^\beta}{m_2^\beta}$ is
$c_1(\beta)u^2(R_1/R_2)^{\beta-1}R_2^{-2}$-spaced, where
$c_1(\beta)>0$ is some positive constant. So from (5.2) we get
\begin{eqnarray*}
{\cal A}_{u,v}(H_1,H_2,R_1, R_2;\delta)&\ll
\frac{H_1H_2}{v^2}\large(1+\frac{\delta R_2^{-\beta}H_1^{-\alpha} }{u^2(R_1/R_2)^{\beta-1}R_2^{-2}}\large)\\
&\ll \frac{H_1H_2}{v^2}+\frac{\delta H_1H_2R_1R_2 }{u^2v^2R_1^\beta
H_1^\alpha   }.
\end{eqnarray*}
Similarly, $\frac{m_2^\beta}{m_1^\beta}$ is
$c_2(\beta)u^2(R_2/R_1)^{\beta-1}R_1^{-2}$-spaced for some positive
constant $c_2(\beta)$, so from (5.3) we get
\begin{eqnarray*}
{\cal A}_{u,v}(H_1,H_2,R_1, R_2;\delta )&\ll
\frac{H_1H_2}{v^2}\large(1+\frac{\delta R_1^{-\beta}H_2^{-\alpha} }{u^2(R_2/R_1)^{\beta-1}R_1^{-2}}\large)\\
&\ll \frac{H_1H_2}{v^2}+\frac{\delta H_1H_2R_1R_2}{u^2v^2R_2^\beta
H_2^\alpha  }.
\end{eqnarray*}

From the above two estimates we get
\begin{eqnarray}
&&\ \ \ \ \ {\cal A}_{u,v}(H_1,H_2,R_1, R_2;\delta)\\&& \ll
\frac{H_1H_2}{v^2}+\frac{\delta
H_1H_2R_1R_2}{u^2v^2}\min(H_1^{-\alpha}R_1^{-\beta} ,
,H_2^{-\alpha}R_2^{-\beta})
\nonumber\\
&&\ll \frac{H_1H_2}{v^2}+\frac{\delta
(H_1H_2)^{1-\alpha/2}(R_1R_2)^{1-\beta} }{u^2v^2},\nonumber
\end{eqnarray}
if we note that $\min(x,y)\leq x^{1/2}y^{1/2}$ for $x>0, y>0.$

Similarly we have
\begin{eqnarray*}
&&  {\cal A}_{u,v}(H_1,H_2,R_1, R_2;\delta) \ll
\frac{R_1R_2}{u^2}+\frac{\delta
(H_1H_2)^{1-\alpha/2}(R_1R_2)^{1-\beta} }{u^2v^2},
\end{eqnarray*}

 which combining (5.4) gives
\begin{eqnarray*}
&&\ \ \ \ \ {\cal A}_{u,v}(D_1,D_2,N_1, N_2;\delta)\\&& \ll
\frac{\delta (H_1H_2)^{1-\alpha/2}(R_1R_2)^{1-\beta} }{u^2v^2}+
\min(\frac{R_1R_2}{u^2},\frac{H_1H_2}{v^2})\\&& \ll \frac{\delta
(H_1H_2)^{1-\alpha/2}(R_1R_2)^{1-\beta}
}{u^2v^2}+\frac{(H_1H_2R_1R_2)^{1/2}}{uv}.
\end{eqnarray*}
Summing over $u$ and $v$ completes the proof of Lemma 5.1.
\end{proof}

\subsection{\bf Estimate of the sum $ S_{a,b}(T)$}\

In this subsection we estimate the sum $S_{a,b}(T).$ For simplicity,
let
$$\eta:=h_1^{\frac{a}{a+b}}r_1^{\frac{b}{a+b}}-
h_2^{\frac{a}{a+b}}r_2^{\frac{b}{a+b}}.$$

 By a splitting argument we get for some
$$(1,1,1,1)\ll(H_1, H_2, R_1, R_2)\ll
(T^{\frac{100(a+b)}{a}}, T^{\frac{100(a+b)}{a}},
T^{\frac{100(a+b)}{b}}, T^{\frac{100(a+b)}{b}})$$ that
\begin{eqnarray}
&& S_{a,b}(T)\ll {\cal L}^4U_{a,b}(T; H_1, H_2, R_1, R_2),
\end{eqnarray}
where
\begin{eqnarray*}
&&U_{a,b}(T; H_1, H_2, R_1, R_2):=
\sum_{3}\frac{1}{(h_1h_2)^{\frac{a+2b}{2a+2b}}(r_1r_2)^{\frac{b+2a}{2a+2b}}}
\min (T^{\frac{1}{a+b}}, |\eta|^{-1} ),\\
&&SC(\Sigma_3): h_j\sim H_j, r_j\sim R_j,(j=1,2), |\eta|<\frac{
1}{10}
h_1^{\frac{a}{2a+2b}}r_1^{\frac{b}{2a+2b}}h_2^{\frac{a}{2a+2b}}r_2^{\frac{b}{2a+2b}}.
\end{eqnarray*}

   By Lemma 5.1 with $(\alpha,\beta)=(a/(a+b),b/(a+b))$  , the contribution of $T^{1/(a+b)}$ is
\begin{eqnarray*}
&& \ll
\frac{T^{\frac{1}{a+b}}}{(H_1H_2)^{\frac{a+2b}{2a+2b}}(R_1R_2)^{\frac{b+2a}{2a+2b}}}
\times {\cal A}(H_1, H_2, R_1,
R_2;T^{-\frac{1}{a+b}})\\
&&\ll 1+ \frac{T^{\frac{1}{a+b}}}{(H_1H_2)^{\frac{
b}{2a+2b}}(R_1R_2)^{\frac{ a}{2a+2b}}}{\cal L}^2.
\end{eqnarray*}

The condition $SC(\Sigma_3)$ implies that $ h_1^ar_1^b\asymp
 h_2^ar_2^b.$ So by the mean value theorem  we get for some $x_0\asymp
 h_1^ar_1^b$ that
\begin{equation}
|\eta|=\frac{1}{a+b}x_0^{\frac{1}{a+b}-1}| h_1^ar_1^b-
h_2^ar_2^b|\gg ( h_1^ar_1^b)^{\frac{1}{a+b}-1},
\end{equation}
which combining the inequality $|\eta|\leq T^{-1/(a+b)}$ gives
$$H_1^aR_1^b\gg T^{1/(a+b-1)}.$$ Similarly we have
$$H_2^aR_2^b\gg
T^{1/(a+b-1)} .$$ If $a\leq b$, then combining the above estimates
  we get that the contribution of $T^{1/(a+b)}$
is
\begin{eqnarray}
&&   \ll 1+ \frac{T^{\frac{1}{a+b}}}{(H_1^aR_1^b)^{\frac{
a}{(2a+2b)b}}(H_2^aR_2^b)^{\frac{ a}{(2a+2b)b}}(H_1H_2)^{\frac{b-
a}{2b}}}{\cal L}^2\nonumber\\ &&\ll
T^{\frac{1}{a+b}-\frac{a}{b(a+b)(a+b-1)}}{\cal L}^2.\nonumber
\end{eqnarray}
If $a>b,$ then the contribution of $T^{1/(a+b)}$ is
\begin{eqnarray}
&&   \ll 1+ \frac{T^{\frac{1}{a+b}}}{(H_1^aR_1^b)^{\frac{
b}{(2a+2b)a}}(H_2^aR_2^b)^{\frac{
b}{(2a+2b)a}}(R_1R_2)^{\frac{a-b}{2a}}}{\cal L}^2\nonumber\\ &&\ll
T^{\frac{1}{a+b}-\frac{b}{a(a+b)(a+b-1)}}{\cal L}^2.\nonumber
\end{eqnarray}
Namely, the the contribution of $T^{1/(a+b)}$ is
\begin{eqnarray}
&&   \ll  T^{\frac{1}{a+b}-\frac{1}{(a+b)(a+b-1)}\times
\frac{\min(a,b)}{\max(a,b)}}{\cal L}^2.
\end{eqnarray}

Now we consider the contribution of $1/|\eta|$. Suppose first that
$a\leq b.$ By a splitting argument the contribution of $1/|\eta|$ is
\begin{eqnarray*}
&& \ll \frac{{\cal L}
}{(H_1H_2)^{\frac{a+2b}{2a+2b}}(R_1R_2)^{\frac{b+2a}{2a+2b}}\delta}
\times {\cal A}(H_1, H_2, R_1,
R_2;\delta)\\
&&\ll  {\cal L} + \frac{1}{(H_1H_2)^{\frac{
b}{2a+2b}}(R_1R_2)^{\frac{ a}{2a+2b}}\delta}{\cal L}^3
\end{eqnarray*}
for some $T^{-1/(a+b)}\ll \delta\ll
H_1^{\frac{a}{2a+2b}}R_1^{\frac{b}{2a+2b}}H_2^{\frac{a}{2a+2b}}R_2^{\frac{b}{2a+2b}}.$
From (5.6) we get
$$\delta\gg ( H_1^aR_1^b)^{\frac{1}{a+b}-1},$$
which combining $T^{-1/(a+b)}\ll \delta$ gives
$$\delta^{-1}\ll \min\left(T^{\frac{1}{a+b}}, ( H_1^aR_1^b)^{\frac{a+b-1}{a+b} }\right).$$
Thus the contribution of $1/|\eta|$ is
\begin{eqnarray*}
&& \ll    {\cal L} + \frac{ 1}{( H_1^aR_1^b)^{\frac{ a}{(2a+2b)b}}(
H_2^aR_2^b)^{\frac{ a}{(2a+2b)b}}(H_1H_2)^{\frac{
b-a}{2b}}\delta}{\cal L}^3\\ &&\ll   {\cal L}  + \frac{1}{(
H_1^aR_1^b)^{\frac{ a}{(a+b)b}}
(H_1H_2)^{\frac{ b-a}{2b}}\delta}{\cal L}^3\nonumber\\
&&\ll    {\cal L} + \frac{ \min\left(T^{\frac{1}{a+b}}, (
H_1^aR_1^b)^{\frac{a+b-1}{a+b} }\right)}{( H_1^aR_1^b)^{\frac{
a}{(a+b)b}} (H_1H_2)^{\frac{ b-a}{2b}}}{\cal L}^3\nonumber\\
&&\ll    {\cal L} +  (H_1H_2)^{-\frac{ b-a}{2b}}\\&&\ \ \ \ \ \ \ \
\ \ \times \min\left(T^{\frac{1}{a+b}}( H_1^aR_1^b)^{-\frac{
a}{(a+b)b}}, ( H_1^aR_1^b)^{\frac{a+b-1}{a+b}-\frac{a}{b(a+b)}
}\right){\cal L}^3\nonumber\\
&&\ll  \left(T^{\frac{1}{a+b}}( H_1^aR_1^b)^{-\frac{
a}{(a+b)b}}\right)^{1-\frac{a}{b(a+b-1)}}\\&&\ \ \ \ \ \ \ \ \
\times \left( ( H_1^aR_1^b)^{\frac{a+b-1}{a+b}-\frac{a}{b(a+b)}
}\right)^{\frac{a}{b(a+b-1)}}{\cal L}^3\nonumber\\&&\ll
  T^{\frac{1}{a+b}-\frac{a}{b(a+b)(a+b-1)}}{\cal L}^3,\nonumber
\end{eqnarray*}
where in the second step we used the fact $ H_1^aR_1^b\asymp
 H_2^aR_2^b.$

Similarly if $a>b,$ we can get that the contribution of $1/|\eta|$
is
\begin{eqnarray}
&&   \ll  T^{\frac{1}{a+b}-\frac{b}{a(a+b)(a+b-1)}}{\cal
L}^3.\nonumber
\end{eqnarray}
Hence combining the above we see that the contribution of $1/|\eta|$
is
\begin{eqnarray}
&&   \ll  T^{\frac{1}{a+b}-\frac{1}{(a+b)(a+b-1)}\times
\frac{\min(a,b)}{\max(a,b)}}{\cal L}^3.
\end{eqnarray}

From (5.5), (5.7) and (5.8) we get the following

{\bf Lemma 5.2.} Suppose $a$ and $ b$ are fixed natural numbers,
then we have
$$S_{a,b}(T)\ll T^{\frac{1}{a+b}-\frac{1}{(a+b)(a+b-1)}\times
\frac{\min(a,b)}{\max(a,b)}}{\cal L}^7
   .$$

\section{\bf Proof of Theorem }

In this section we shall prove our Theorem . It suffices to evaluate
the integral $\int_T^{2T}\Delta^2_{a,b}(x)dx$ for $T\geq 10.$ We
always suppose $a<b.$

\subsection{\bf Mean square of $\Delta^{*}_{a,b}(x,z)$}\

Suppose $H$ and $z$ are parameters such that $$T^\varepsilon\ll H\ll
T^{100(a+b)}, \ \  T^\varepsilon\ll z\leq \min(H^{a+b},T^{a/b}){\cal
L}^{-b-b^2/a-1}.$$ In this subsection we study the mean square of
 $\Delta^{*}_{a,b}(x,z).$

  By the elementary formula
  \begin{equation}
\cos u\cos v=\frac{1}{2}(\cos {(u-v)}+\cos {(u+v)})
 \end{equation}
we may write
\begin{eqnarray}
&&\ \ \ \ \ \ |\Delta^{*}_{a,b}(x,z) |^2\\&&=
\frac{c_1^2(a,b)}{\pi^2}x^{\frac{1}{(a+b)}}\sum_{1\leq n\leq z
}\sum_{1\leq m\leq z }g_{a,b}(n)g_{a,b}(m)\nonumber\\&&\ \ \ \ \ \ \
\times\cos\left(2\pi
c_2(a,b)x^{\frac{1}{a+b}}n^{\frac{1}{a+b}}-\frac{\pi}{4}\right)
\cos\left(2\pi
c_2(a,b)x^{\frac{1}{a+b}}m^{\frac{1}{a+b}}-\frac{\pi}{4}\right)\nonumber\\
&&=S_1(x)+S_2(x)+S_3(x),\nonumber
\end{eqnarray}
where
\begin{eqnarray*}
S_1(x)&&=\frac{c_1^2(a,b)}{2\pi^2}x^{\frac{1}{(a+b)}}\sum_{1\leq
n\leq z
}g^2_{a,b}(n), \\
S_2(x)&&=\frac{c_1^2(a,b)}{2\pi^2}x^{\frac{1}{(a+b)}}\sum_{\stackrel{n\leq z,m\leq z}{n\not= m}}g_{a,b}(n)g_{a,b}(m) \\
&&\ \ \ \ \ \ \ \ \times \cos\left(2\pi
c_2(a,b)x^{\frac{1}{a+b}}(n^{\frac{1}{a+b}}-m^{\frac{1}{a+b}}
)\right)\\
 S_3(x)&&=\frac{c_1^2(a,b)}{2\pi^2}x^{\frac{1}{(a+b)}}\sum_{ n\leq z,m\leq z }g_{a,b}(n)g_{a,b}(m) \\
&&\ \ \ \ \ \ \ \ \times \sin\left(2\pi
c_2(a,b)x^{\frac{1}{a+b}}(n^{\frac{1}{a+b}}+m^{\frac{1}{a+b}}
)\right) .
\end{eqnarray*}

 By (4.6) we get
 \begin{equation}
\int_T^{2T}S_1(x)dx=\frac{c_1^2(a,b)}{2\pi^2} \sum_{n=1}^\infty
g_{a,b}^2(n)\int_T^{2T}x^{\frac{1}{a+b}}dx+O(T^{\frac{1+a+b}{a+b}}z^{-\frac{a}{(a+b)b}}).
\end{equation}

By the first derivative test we get
\begin{eqnarray*}
\int_T^{2T}S_3(x)dx&&\ll \sum_{ n\leq z,m\leq z
}g_{a,b}(n)g_{a,b}(m)
\frac{T}{n^{\frac{1}{a+b}}+m^{\frac{1}{a+b}}}\\
&&\ll T\sum_{ n\leq z,m\leq z
}g_{a,b}(n)g_{a,b}(m)\frac{1}{n^{\frac{1}{2(a+b)}}m^{\frac{1}{2(a+b)}}} \\
&&\ll T\left(\sum_{ n\leq z
}\frac{g_{a,b}(n)}{n^{\frac{1}{2(a+b)}}}\right)^2,
\end{eqnarray*}
where in the second step we used the well-known inequality
$\alpha^2+\beta^2\geq 2\alpha\beta.$ By Euler's product we have for
$\Re s>1$ that
\begin{eqnarray*}
\sum_{n=1}^\infty\frac{g_{a,b}(n)}{n^s}=\zeta(as+\frac{a+2b}{2a+2b})\zeta(bs+\frac{2a+b}{2a+2b}),
\end{eqnarray*}
which can be continued meromorphically to the whole complex plane
and has a double pole at $s=1/2(a+b).$ Thus Perron's formula implies
that
\begin{eqnarray*}
\sum_{n\leq X} g_{a,b}(n)\ll X^{1/2(a+b)}\log X
\end{eqnarray*}
for any $X>2,$ namely,
\begin{eqnarray}
\sum_{n\leq X} g_{a,b}(n)n^{-1/2(a+b)}\ll \log^2 X.
\end{eqnarray}
From the above we get
\begin{eqnarray}
\int_T^{2T}S_3(x)dx\ll  T{\cal L}^4 .
\end{eqnarray}

Now we consider the contribution of $S_2(x).$ Write
\begin{eqnarray}
S_2(x)=S_{21}(x)+S_{22}(x),
\end{eqnarray}
where
\begin{eqnarray*}
S_{2j}(x)&&=\frac{c_1^2(a,b)}{2\pi^2}x^{\frac{1}{(a+b)}}\sum_{2j}g_{a,b}(n)g_{a,b}(m) \\
&&\ \ \ \ \ \ \ \ \times \cos\left(2\pi
c_2(a,b)x^{\frac{1}{a+b}}(n^{\frac{1}{a+b}}-m^{\frac{1}{a+b}}
)\right)\ \  (j=1,2),\\
SC(\Sigma_{21}):&& n,m\leq z,
|n^{\frac{1}{a+b}}-n^{\frac{1}{a+b}}|\geq
\frac{(nm)^{1/(2a+2b)}}{10},\\
SC(\Sigma_{22}):&& n,m\leq z, 0<
|n^{\frac{1}{a+b}}-n^{\frac{1}{a+b}}|< \frac{(nm)^{1/(2a+2b)}}{10}.
\end{eqnarray*}

Similar to the case $S_3(x),$ we have
\begin{eqnarray*}
\int_T^{2t}S_{21}(x)dx\ll T{\cal L}^4.
\end{eqnarray*}

By the first derivative test  and Lemma 5.2 we get
\begin{eqnarray*}
\int_T^{2T}S_{22}(x)dx&&\ll
T\sum_{22}g_{a,b}(n)g_{a,b}(m)\min\left(T^{\frac{1}{a+b}},
\frac{1}{|n^{\frac{1}{a+b}}-n^{\frac{1}{a+b}}|}\right)\\
&&\ll TS_{a,b}(T)\ll
T^{1+\frac{1}{a+b}-\frac{a}{b(a+b)(a+b-1)}}{\cal L}^7.
\end{eqnarray*}

From the above two estimates we get
\begin{eqnarray}
\int_T^{2T}S_{2}(x)dx  \ll
T^{1+\frac{1}{a+b}-\frac{a}{b(a+b)(a+b-1)}}{\cal L}^7.
\end{eqnarray}

From (6.2), (6.3), (6.5) and (6.7) we have
 \begin{eqnarray}
&&\int_T^{2T}|\Delta^{*}_{a,b}(x,z)|^2dx=\frac{c_1^2(a,b)}{2\pi^2}
\sum_{n=1}^\infty
g^2_{a,b}(n)\int_T^{2T}x^{\frac{1}{a+b}}dx\\
&&\ \ \ \ \ \ \ \ \ \ \ \ \ \ \ \ \ \ \ \ \ \ \ \ \ \ \ \
+O(T^{\frac{1+a+b}{a+b}}z^{-\frac{a}{(a+b)b}}+T^{\frac{1+a+b}{a+b}-\frac{a}{b(a+b)(a+b-1)}}{\cal
L}^7 ).\nonumber
\end{eqnarray}

\subsection{\bf Mean squares of $R_{12}^{*}(a,b;x)$ and  $R_{12}^{*}(b,a;x)$}\

In this subsection we shall study the mean squares of
$R_{12}^{*}(a,b;x)$ and  $R_{12}^{*}(b,a;x).$

Recall that
\begin{eqnarray*}
R_{12}^{*}(a,b; x)&&=
\frac{c_1(a,b)}{\pi}x^{\frac{1}{2(a+b)}}\sum_{z\leq n\leq
H^an_{J+1,H}(a,b)}g(a,b;n,H,J)\\&&\ \ \ \ \ \ \ \times\cos\left(2\pi
c_2(a,b)x^{\frac{1}{a+b}}n^{\frac{1}{a+b}}-\frac{\pi}{4}\right).
\end{eqnarray*}
 Thus
\begin{eqnarray*}
&&\ \ \ \ \ \ \ |R_{12}^{*}(a,b; x)|^2\\&&\ll
x^{\frac{1}{a+b}}\sum_{z\leq n,m\leq H^an_{J+1,H}(a,b)}g(a,b;n,H,J)
g(a,b;m,H,J)\\&&\ \ \ \ \ \  \ \ \ \ \ \ \ \ \ \ \ \ \ \ \ \ \
\times e\left(
c_2(a,b)x^{\frac{1}{a+b}}(n^{\frac{1}{a+b}}-m^{\frac{1}{a+b}})
\right)\\
&&=S_4(x)+S_5(x),
\end{eqnarray*}
where
\begin{eqnarray*}
&&S_4(x) =x^{\frac{1}{a+b}}\sum_{z\leq n\leq
H^an_{J+1,H}(a,b)}g^2(a,b;n,H,J),\\
&&S_5(x)=x^{\frac{1}{a+b}}\sum_{\stackrel{z\leq n,m\leq
H^an_{J+1,H}(a,b)}{n\not= m}}g(a,b;n,H,J) g(a,b;m,H,J)\\&&\ \ \ \ \
\ \ \ \ \ \ \ \ \ \ \ \ \ \ \ \ \ \ \times e\left(
c_2(a,b)x^{\frac{1}{a+b}}(n^{\frac{1}{a+b}}-m^{\frac{1}{a+b}})
\right)
\end{eqnarray*}

It is easy to see that
\begin{equation}
g(a,b;n,H,J)\leq g(a,b;n)\leq g_{a,b}(n).
\end{equation}
By (4.7) we get
\begin{eqnarray}
\int_T^{2T}S_4(x)dx\ll T^{1+\frac{1}{a+b}}z^{-\frac{a}{b(a+b)}}
\end{eqnarray}.

By the first derivative test  and  (6.9)  we get
\begin{eqnarray}
\int_T^{2T}S_5(x)dx&&\ll T\sum_{\stackrel{z\leq n,m\leq H^aT}{n\not=
m}}g_{a,b}(n) g_{a,b}(m)\\&&\ \ \ \ \ \ \ \times
\min\left(T^{\frac{1}{a+b}},\frac{1}{|n^{\frac{1}{a+b}}-m^{\frac{1}{a+b}}|}\right)\nonumber\\
&&\ll T(\Sigma_4+\Sigma_5),\nonumber
 \end{eqnarray}
  where
\begin{eqnarray*}
&&\Sigma_4:=\sum_{4}g_{a,b}(n) g_{a,b}(m)
\min\left(T^{\frac{1}{a+b}},\frac{1}{|n^{\frac{1}{a+b}}-m^{\frac{1}{a+b}}|}\right),\\
&&\Sigma_5:=\sum_{5}g_{a,b}(n) g_{a,b}(m)
\min\left(T^{\frac{1}{a+b}},\frac{1}{|n^{\frac{1}{a+b}}-m^{\frac{1}{a+b}}|}\right),\\
&&SC(\Sigma_4): z<n,m\leq H^aT,  0<
|n^{\frac{1}{a+b}}-m^{\frac{1}{a+b}}|\leq
\frac{1}{10}n^{\frac{1}{2a+2b}}m^{\frac{1}{2a+2b}},\\
&&SC(\Sigma_5): z<n,m\leq H^aT,
|n^{\frac{1}{a+b}}-m^{\frac{1}{a+b}}|\geq
\frac{1}{10}n^{\frac{1}{2a+2b}}m^{\frac{1}{2a+2b}}.
 \end{eqnarray*}

By Lemma 5.2 we have
\begin{equation}
\Sigma_4\ll S_{a,b}(T)\ll
T^{\frac{1}{a+b}-\frac{a}{b(a+b)(a+b-1)}}{\cal L}^7.
\end{equation}

For $\Sigma_5$ by (6.4) we get
\begin{equation}
\Sigma_5 \ll  (\sum_{n\leq H^aT}g_{a,b}(n)n^{-1/(2a+2b)})^2\ll
\log^4 X.
\end{equation}

From (6.10)-(6.13) we have
\begin{equation}
\int_T^{2T}|R_{12}^{*}(a,b; x)|^2dx\ll
T^{1+\frac{1}{a+b}}z^{-\frac{a}{b(a+b)}}+
T^{\frac{1+a+b}{a+b}-\frac{a}{b(a+b)(a+b-1)}}{\cal L}^7.
\end{equation}

Similarly, we have
\begin{equation}
\int_T^{2T}|R_{12}^{*}(b,a; x)|^2dx\ll
T^{1+\frac{1}{a+b}}z^{-\frac{a}{b(a+b)}}+
T^{\frac{1+a+b}{a+b}-\frac{a}{b(a+b)(a+b-1)}}{\cal L}^7.
\end{equation}

\subsection{\bf Mean squares of $R_2(a,b;x)$ and $R_2(b,a;x)$}\

In this subsection we shall study the mean squares of $R_2(a,b;x)$
and $R_2(b,a;x)$. Recall that
\begin{eqnarray*}
 R_2(a,b; x)=O\left(\sum_{m\leq
x^{1/(a+b)}}\min\left(1,\frac{1}{H\Vert \frac{x^{1/a}}{m^{b/a}}
\Vert}\right) \right).
\end{eqnarray*}
So we have
\begin{eqnarray*}
&&\ \ \ \ \ \ \ \ \ \ \int_T^{2T}R_2(a,b; x)dx\\&&\ll \sum_{m\leq
{2T}^{1/(a+b)}}\int_T^{2T}\min\left(1,\frac{1}{H\Vert
\frac{x^{1/a}}{m^{b/a}} \Vert}\right)dx\nonumber\\
&&\ll \sum_{m\leq
{2T}^{1/(a+b)}}\int_{T^{1/a}/m^{b/a}}^{(2T)^{1/a}/m^{b/a}}\min\left(1,\frac{1}{H\Vert
 u\Vert}\right)m^bu^{a-1}du\nonumber\\
 &&\ll \sum_{m\leq
{2T}^{1/(a+b)}}m^b\left(\frac{T^{1/a}}{m^{b/a}}\right)^{a-1}\int_{T^{1/a}/m^{b/a}}^{(2T)^{1/a}/m^{b/a}}\min\left(1,\frac{1}{H\Vert
 u\Vert}\right) du\nonumber\\
 &&\ll \sum_{m\leq
{2T}^{1/(a+b)}}m^b\left(\frac{T^{1/a}}{m^{b/a}}\right)^{a-1}\frac{T^{1/a}}{m^{b/a}}\int_{0}^{1}\min\left(1,\frac{1}{H\Vert
 u\Vert}\right) du\nonumber\\
 &&\ll T \sum_{m\leq
{2T}^{1/(a+b)}}\int_{0}^{1}\min\left(1,\frac{1}{H\Vert
 u\Vert}\right) du\nonumber\\
 &&\ll T^{1+1/(a+b)}\int_{0}^{1/2}\min\left(1,\frac{1}{H\Vert
 u\Vert}\right) du\nonumber\\
 &&\ll T^{1+1/(a+b)}\left(\int_0^{1/H}du+\int_{1/H}^{1/2}\frac{1}{Hu}du\right)\\
 &&\ll T^{1+1/(a+b)}H^{-1}{\cal L},\nonumber
\end{eqnarray*}
which combining the trivial bound $R_2(a,b;x)\ll T^{1/(a+b)}$ gives
\begin{eqnarray}
 \int_T^{2T}R_2^2(a,b; x)dx\ll T^{1+2/(a+b)}H^{-1}{\cal L}
\end{eqnarray}

Similarly we have
\begin{eqnarray}
&&  \int_T^{2T}R_2^2(b,a; x)dx\ll T^{1+2/(a+b)}H^{-1}{\cal L}.
\end{eqnarray}

\subsection{\bf Completion of the proof of Theorem }\

In this subsection we complete the proof of Theorem . We take
$H=T^{10(a+b)}$ and  $z=T^{a/b}{\cal L}^{-b-b^2/a-1}.$
 From (3.11), (6.11)-(6.14) we get
\begin{eqnarray}
\int_T^{2T}E^2_{a,b}(x)dx\ll
T^{\frac{1+a+b}{a+b}-\frac{a}{b(a+b)(a+b-1)}}{\cal L}^7,
\end{eqnarray}
which combining (6.8) and Cauchy's inequality gives
\begin{equation}
\int_T^{2T}E_{a,b}(x)\Delta^{*}_{a,b}(x,z)dx\ll
T^{\frac{1+a+b}{a+b}-\frac{a}{2b(a+b)(a+b-1)}}{\cal L}^{7/2}.
\end{equation}

From (6.8), (6.18) and  (6.19) we get
\begin{eqnarray}
&&\int_T^{2T}\Delta^2_{a,b}(x)dx=\frac{c_1^2(a,b)}{2\pi^2}
\sum_{n=1}^\infty
g^2_{a,b}(n)\int_T^{2T}x^{\frac{1}{a+b}}dx\\
&&\ \ \ \ \ \ \ \ \ \ \ \ \ \ \ \ \ \ \ \ \ \ \ \ \ \ \ \ +O(
T^{\frac{1+a+b}{a+b}-\frac{a}{2b(a+b)(a+b-1)}}{\cal L}^{7/2}
).\nonumber
\end{eqnarray}

From (6.20) and a splitting argument we get
\begin{eqnarray}
\int_1^{T}\Delta^2_{a,b}(x)dx&&=\frac{c_1^2(a,b)}{2\pi^2}
\sum_{n=1}^\infty
g^2_{a,b}(n)\int_1^{T}x^{\frac{1}{a+b}}dx\\
&&\ \ \ \ \ \ \   +O(
T^{\frac{1+a+b}{a+b}-\frac{a}{2b(a+b)(a+b-1)}}{\cal L}^{7/2}
)\nonumber\\
&&=\frac{c_1^2(a,b)(a+b)}{2(1+a+b)\pi^2} \sum_{n=1}^\infty
g^2_{a,b}(n) T^{\frac{1+a+b}{a+b}}\nonumber\\
&&\ \ \ \ \ \ \   +O(
T^{\frac{1+a+b}{a+b}-\frac{a}{2b(a+b)(a+b-1)}}{\cal L}^{7/2}
)\nonumber\\&& =c_{a,b}T^{\frac{1+a+b}{a+b}} +O(
T^{\frac{1+a+b}{a+b}-\frac{a}{2b(a+b)(a+b-1)}}{\cal
L}^{7/2}).\nonumber
\end{eqnarray}

Authors' addresses:

Wenguang Zhai,\\
School of Mathematical Sciences,\\
Shandong Normal University,\\Jinan,  Shandong, 250014,\\
P.R.China\\
 E-mail:  zhaiwg@hotmail.com

Xiaodong Cao\\ Dept. of Mathematics and Physics,\\ Beijing Institute
of
Petro-Chemical Technology,\\
Beijing, 102617, P.R. China\\
Email: caoxiaodong@bipt.edu.cn


\begin{thebibliography}{99}

\bibitem{Cx}X. D. Cao, On the general divisor problem(in Chinese),
Acta Mathematica sinica, Vol.{\bf 36}(5)(1993), 644-653.

\bibitem{CK}K. Corr\'adi and I. Katai, A comment on K. S. Ganggadharan's paper
entitled  `` Two classical lattice point problems'', Magyar Tud.
Akad. Mat. Fiz. Oszt. K\"ozl . {\bf 17}(1967), 89-97.

\bibitem{Cr}H. Cram\'er, \"Uber zwei S\"atze von Herrn G. H. Hardy,
Math. Z.{\bf 15}(1922), 201-210.

\bibitem{FI} E. Fouvry and H. Iwaniec, Exponential sums with monomials,
 J. of Number Theory, Vol. {\bf  33}(3)(1989), 311---333 .

\bibitem{Ha1}J. L. Hafner, New omega theorems for two classical lattice point pro
blems, Invent. Math. {\bf 63}(1981), 181-186.

\bibitem{Ha2}J. L. Hafner, New omega
results in a weighted divisor problem, J. of Number Theory Vol. {\bf
28}(1988), 240-257.

\bibitem{H1}D. R. Heath-Brown, The distribution and moments of the error term
in the Dirichlet divisor problem, Acta Arith. {\bf 60} (1992),
389-415.

\bibitem{H2}D. R. Heath-Brown, The Piatetski-Shapiro prime theorem, J. of Number
theory, Vol. {\bf 16}(1983), 242-266.

\bibitem{Hu}M. N. Huxley,  Exponential sums and
Lattice points III, Proc. London Math. Soc., Vol. {\bf 87}(3)(2003),
591-609.

\bibitem{I1}A. Ivi\'c, Lectures on mean values of the Riemann zeta-function,
Lectures On Math. and Physics 82, Tata Inst. Fund. Res., Bombay, 1991.



\bibitem{I3}A. Ivi\'c, The Riemann zeta-function. John, Wiley and Sons,
1985.

\bibitem{I4} A. Ivi\'c, The general divisor problem, J. of Number
Theory {\bf 27}(1987), 73-91.

\bibitem{IS}A. Ivi\'c  and P. Sargos, On the higher power moments
of the error term in the divisor problem, Illinois Journal of Math.
{\bf 51}(2)(2007), 353-377.

\bibitem{Kr1}E. Kr\"{a}tzel, Ein Teilerprolem, J. Reine Angew. Math.
{\bf 235}(1969), 150-174.

\bibitem{Kr2}  E. Kr\"atzel, Lattice points, Deutsch. Verlag Wiss., Berlin, 1988.

\bibitem{KN}M. K\"{u}hleitner and W. G. Nowak, The asymptotic behaviour of the mean-square of fractional
part sums, Proc. Edinb. Math. Soc. 43 (2000), 309-323.

\bibitem{Mi}S. H. Min, The methods of number theory(in Chinese),
Science Press, Beijing: 1981.


\bibitem{P}E. Preissmann, Sur la moyenne quadratique du terme de reste du probl¨¦me du cercle,
C. R. Acad. Sci. Paris S¨¦r. I 306 (1988), 151-154.

\bibitem{Ra}R.A. Rankin, Van Der Corput's method and the theory of
exponent pairs, Quart. J. Math. Oxford (2){\bf 6}(1955), 147-153.

\bibitem{Ri}H.-E. Richert, \"{U}ber die Anzahl Gruppen gegebener
Ordnung (I), Math. Z. {\bf 56}(1952), 21-32.

\bibitem{S}P. G. Schmidt, Absch\"{a}tzungen bei unsymmetrischen
Gitterpunktproblem, Dissertation , G\"{o}ttingen, 1964.


\bibitem{To}K. C. Tong, On divisor problem III, Acta math. Sinica {\bf 6}
(1956), 515-541.

\bibitem{Ts}Kai-Man Tsang, Higher-power moments of $\Delta(x), E(t)$ and $P(x)$,
Proc. London Math. Soc.(3){\bf 65}(1992), 65-84.

\bibitem{V} I. M. Vinogradov, Special variants of the method of trigonometric sums,
 (Nauka, Moscow), 1976; English transl. in his Selected works (Springer-Verlag),
 1985.

\bibitem{Zw1}Wenguang Zhai , On higher-power moments of $\Delta(x),$ Acta
Arith. Vol.{\bf 112}(2004), 1-24.


\bibitem{Zw2} Wenguang Zhai, On higher-power moments of $\Delta(x) $(II), Acta Arith.
Vol.{\bf 114}(2004), 35-54.

\bibitem{Zw3} Wenguang Zhai, On higher-power moments of $\Delta(x)
$(III), Acta Arith. {\bf 118} (2005), 263--281.

\end{thebibliography}
\end{document}